\documentclass{amsart}
\usepackage{amsmath,amssymb,cite,geometry,bm}
\usepackage[colorlinks,citecolor=blue,linkcolor=blue,urlcolor=blue]{hyperref}
\newtheorem{theorem}{Theorem}[section]

\newtheorem{remark}[theorem]{Remark}
\numberwithin{equation}{section}
\allowdisplaybreaks[4]

\begin{document}

\title{Two-term spectral asymptotics in linear elasticity on a Riemannian manifold}

\author{Genqian Liu}
\address{School of Mathematics and Statistics, Beijing Institute of Technology, Beijing 100081, China}
\email{liugqz@bit.edu.cn}


\begin{abstract}
\ In this note, by explaining two key methods that were employed in \cite{Liu-21} and by giving some remarks, we show that the proof of Theorem 1.1 in \cite{Liu-21} is a rigorous  proof based on theory of strongly continuous semigroups and pseudodifferential operators. All remarks and comments to paper \cite{Liu-21}, which were given by Matteo Capoferri, Leonid Friedlander, Michael Levitin and Dmitri Vassiliev in \cite{CaFrLeVa-22}, are incorrect. The so-called ``numerical counter-examples'' in \cite{CaFrLeVa-22} are useless examples for the two-term asymptotics of the counting functions of the elastic eigenvalues. Clearly, the conclusion and the proof of \cite{Liu-21} are completely correct.
\end{abstract}

\maketitle
\vskip 0.19 true cm
\section{Two key methods}

 \vskip 0.39 true cm

Let $(\Omega,g)$ be a compact smooth Riemannian $n$-manifold with smooth boundary $\partial \Omega$. Let $P_g$ be the Navier--Lam\'{e} operator (see \cite{Liu-19} and \cite{Liu-21}):
\begin{eqnarray} \label {1-1} P_g\mathbf{u}:=\mu \nabla^* \nabla \mathbf{u} -(\mu +\lambda) \,\mbox{grad}\; \mbox{div}\, \mathbf{u} -\mu \, \mbox{Ric} (\mathbf{u}), \;  \;\; \mathbf{u}=(u^1, \cdots, u^n),\end{eqnarray}
where $\mu$ and $\lambda$ are Lam\'{e} parameters satisfying $\mu>0$ and $\mu+\lambda\ge 0$,  $\nabla^* \nabla $ is the Bochner Laplacian (see (2.11) of \cite{Liu-21}), $\mbox{div}$ and $\mbox{grad}$ are the usual divergence and gradient operators, and  \begin{eqnarray} \label{18/12/22} \mbox{Ric} (\mathbf{u})= \big(\sum\limits_{l=1}^nR^{\,1}_{l} u^l,  \sum\limits_{l=1}^n R^{\,\,2}_{\,l} u^l, \cdots, \sum\limits_{l=1}^n R^{\,\,n}_{\,l} u^l\big)\end{eqnarray} denotes the action of Ricci tensor $\mbox{R}_l^{\;j}:=\sum_{k=1}^n R^{k\,\,j}_{\,lk}$ on $\mathbf{u}$.

We denote by $P_g^-$ and $P_g^+$ the Navier--Lam\'{e} operators with zero Dirichlet and Neumann boundary conditions, respectively. The Dirichlet boundary condition is $\mathbf{u}\big|_{\partial \Omega}$, and the Neumann boundary condition is \begin{eqnarray*}\frac{\partial \mathbf{u}}{\partial \boldsymbol{\nu}}:=2\mu  (\mbox{Def}\, \mathbf{u})^\# \boldsymbol{\nu} + \lambda (\mbox{div}\, \mathbf{u} )\boldsymbol{\nu}\;\;\mbox{on}\;\; \partial \Omega,\end{eqnarray*} where $\mbox{Def}\, \mathbf{u}= \frac{1}{2} (\nabla \mathbf{u} +\nabla \mathbf{u}^T)$,  $\,\#$ is the sharp operator (for a tensor) by raising index, and $\boldsymbol{\nu}$ is the unit outer normal to $\partial \Omega$.

Since $P_g^-$ (respectively, $P_g^+$) is an unbounded, self-adjoint and positive (respectively,  nonnegative) operator in $[H^1_0(\Omega)]^n$ (respectively, $[H^1(\Omega)]^n$) with discrete spectrum $0< \tau_1^- < \tau_2^- \le \cdots \le \tau_k^- \le \cdots \to +\infty$ (respectively, $0\le \tau_1^+ < \tau_2^+ \le \cdots \le \tau_k^+ \le \cdots \to +\infty$), one has
\begin{eqnarray} \label{1-4} P_g^\mp {\mathbf{u}}_k^\mp =\tau_k^\mp {\mathbf{u}}_k^\mp,\end{eqnarray} where ${\mathbf{u}}_k^-\in  [H^1_0(\Omega)]^n$ (respectively, ${\mathbf{u}}_k^+\in  [H^1(\Omega)]^n$) is the eigenvector corresponding to eigenvalue $\tau_k^{-}$ (respectively, $\tau_k^{+}$).

\vskip 0.10 true cm
In \cite{Liu-21}, we obtained the following result:

\vskip 0.22 true cm

\begin{theorem}
Let $(\Omega,g)$ be a compact smooth Riemannian manifold of dimension $n$ with smooth boundary $\partial \Omega$, and let $0< \tau_1^-< \tau_2^- \le \tau^-_3\le \cdots \le \tau_k^- \le \cdots$ (respectively, $0\le \tau_1^+ < \tau_2^+ \le \tau_3^+ \le \cdots \le \tau_k^+ \le \cdots $) be the eigenvalues of the Navier--Lam\'{e} operator $P_g^-$ (respectively, $P_g^+$) with respect to the zero Dirichlet (respectively, Neumann) boundary condition. Then
\begin{eqnarray} \label{1-7} &&  \sum_{k=1}^\infty e^{-\tau_k^\mp t}  = \mbox{Tr}\,(e^{-tP_g^\mp})=\bigg[ \frac{n-1}{(4\pi \mu t)^{n/2}}
 + \frac{1}{(4\pi (2\mu+\lambda) t)^{n/2}}\bigg] {\mbox{Vol}}\,(\Omega) \\
&&  \quad \, \mp \frac{1}{4} \bigg[  \frac{n-1}{(4\pi \mu t)^{(n-1)/2}}
 +  \frac{1}{(4\pi (2\mu+\lambda) t)^{(n-1)/2}}\bigg]{\mbox{Vol}}\,(\partial\Omega)+O(t^{{1-n}/2})\quad\;\; \mbox{as}\;\; t\to 0^+.\nonumber\end{eqnarray}
 Here ${\mbox{Vol}}\,(\Omega)$ denotes the $n$-dimensional volume of $\Omega$,  ${\mbox{Vol}}\, (\partial\Omega)$ denotes the $(n-1)$-dimensional volume of $\partial \Omega$.
\end{theorem}

 \vskip 0.2 true cm

According to the proof of Theorem 1.1, in particular, by the expressions of $\mathbf{q}_{-3}$ and $\mathbf{q}_{-4}$
in (3.10) of \cite{Liu-21}, we actually proved the following result (see, ii) of  Remark 4.1 in \cite{Liu-21}):
\begin{theorem}
  Let $(\Omega,g)$ be a compact smooth Riemannian manifold of dimension $n$ with smooth boundary $\partial \Omega$. Then
 \begin{eqnarray} \label{2022.7.30-1} &&  \sum_{k=1}^\infty e^{-\tau_k^\mp t} = t^{-n/2} \big[a_0+a_1^\mp t^{1/2} + a_2^\mp  t +\cdots +a_m^\mp t^{m/2} +O(t^{(m+1)/2})\big]\;\,\mbox{as}\;\, t\to 0^+,\nonumber\end{eqnarray}
where all coefficients
$a_0$, $\,a_k^{\mp}$, ($k=1,2,3,\cdots)$ can be explicitly obtained by an algorithm  (see the proof of \cite{Liu-21}). In particular,
\begin{eqnarray} && a_0= \bigg[\frac{n-1}{(4\pi \mu )^{n/2}}
 + \frac{1}{(4\pi (2\mu+\lambda) )^{n/2}}\bigg] {\mbox{Vol}}\,(\Omega),\\
 && a_1^\mp=\mp \frac{1}{4} \bigg[  \frac{n-1}{(4\pi \mu )^{(n-1)/2}}
 +  \frac{1}{(4\pi (2\mu+\lambda) )^{(n-1)/2}}\bigg]{\mbox{Vol}}\,(\partial\Omega),\\
 && \label{2022.7.31-6}a_2^\mp= \bigg[  \frac{n-1}{(4\pi \mu )^{n/2}}
 +  \frac{1}{(4\pi (2\mu+\lambda) )^{n/2}}\bigg] \left( \frac{1}{3} \int_\Omega R \,dV -\frac{1}{6} \int_{\partial \Omega} H \, ds\right),\end{eqnarray}
where $R$ is the scalar curvature of $\Omega$ and $H$ is the mean curvature of $\partial \Omega$. \end{theorem}

 \vskip 0.38 true cm

In the proof of Theorem 1.1 in \cite{Liu-21}, we have used two key methods: the ``method of image'' and the ``explicitly calculable method'' for the full  symbol of the resolvent operator $(\tau I-P_g)^{-1}$. In this note, by explaining the two key methods again, and by giving some remarks, we will show that our conclusion and the proof of Theorem 1.1 in \cite{Liu-21} are correct.

\vskip 0.25 true cm

In order to explain the ``method of image'', we first introduce its original idea (see \cite{Liu-21.1} or \cite{Liu-21}) for calculating the Green function of a parabolic  (elastic) system with zero Dirichlet (or Neumann) boundary condition in the upper half-space $\mathbb{R}^{n}_{+}:=\{x=(x',x_n)\in \mathbb{R}^n \big| x_n>0\}$, where $x^{\prime}=(x_1,\dots,x_{n-1})$.

In the upper half-space $\mathbb{R}^n_{+}$, the Lam\'{e} operator $P$ can be written as  \begin{eqnarray*}  P\mathbf{u}\!\!&=&\!\!-\mu\Delta \mathbf{u}-(\lambda+\mu)\nabla(\nabla\cdot \mathbf{u})\\
\!\!\!& =&\!\!\!-\mu \begin{pmatrix} \Delta & \cdots & 0 & 0\\
\vdots & \ddots & \vdots & \vdots\\
0& \cdots & \Delta & 0\\
0 & \cdots & 0 & \Delta\end{pmatrix}\begin{pmatrix} u_1\\
\vdots \\ u_{n-1} \\  u_n\end{pmatrix} -(\lambda+\mu) \begin{pmatrix} \partial^2_{x_1 x_1} & \cdots & \partial^2_{x_1x_{n-1}} & \partial^2_{x_1x_{n}} \\
\cdots& \cdots &\cdots& \cdots\\
\partial^2_{x_{n-1} x_1} & \cdots & \partial^2_{x_{n-1}x_{n-1}} & \partial^2_{x_{n-1}x_{n}}\\
\partial^2_{x_n x_1} & \cdots & \partial^2_{x_nx_{n-1}} & \partial^2_{x_n x_{n}}  \end{pmatrix}\begin{pmatrix} u_1\\
 \vdots \\ u_{n-1} \\ u_n\end{pmatrix}, \end{eqnarray*}
where $\Delta=\sum_{j=1}^n \partial^2_{x_j x_j}$ is the standard Laplacian. Let  $\varsigma: (x_1, \cdots, x_{n-1}, x_n)\mapsto (x_1, \cdots, x_{n-1}, -x_n)$ be the reflection with respect to the $x'$-hyperplane  (here we always assume $x_n\ge 0$).  Then by $\varsigma$,  from the given upper half-space $\mathbb{R}^n_+$ we can get the lower half-space $\mathbb{R}^n_{-}$ (and hence  getting  $\mathbb{\bar R}^n_-$ from $\mathbb{\bar R}^n_+$). In addition, by this reflection $\varsigma$,
the differential operators $\partial_{x_1}$, $\cdots$, $\partial_{x_{n-1}}$, $\partial_{x_n}$ (defined on the upper half-space) are changed into $\partial_{x_1}$, $\cdots$, $\partial_{x_{n-1}}$, $-\partial_{x_n}$ (defined on the lower half-space), respectively.
 We define \begin{eqnarray} \label{2022.10.18-2} \mathcal{P}=\left\{\begin{array}{ll} \! P \;\;\; \;\;\;\,\mbox{on the upper half-space} \;\, \mathbb{\bar R}^n_+\\
 \!P^- \;\;\;\; \mbox{on the lower half-space} \;\, \mathbb{R}^n_-, \end{array} \right.\end{eqnarray}
where \begin{eqnarray*} \label{2022.10.6} P^-:=- \mu \begin{pmatrix} \Delta & \cdots & 0 & 0\\
\vdots & \ddots & \vdots & \vdots\\
0& \cdots & \Delta & 0\\
0 & \cdots & 0 & \Delta\end{pmatrix}\begin{pmatrix} u_1\\
 \vdots \\ u_{n-1} \\ u_n\end{pmatrix} -(\lambda+\mu) \begin{pmatrix} \partial^2_{x_1 x_1} & \cdots & \partial^2_{x_1x_{n-1}} & -\partial^2_{x_1x_{n}} \\
\vdots& \ddots &\vdots& \vdots\\
\partial^2_{x_{n-1} x_1} & \cdots & \partial^2_{x_{n-1}x_{n-1}} & -\partial^2_{x_{n-1}x_{n}}\\
-\partial^2_{x_n x_1} & \cdots & -\partial^2_{x_nx_{n-1}} & \partial^2_{x_n x_{n}}  \end{pmatrix}\begin{pmatrix} u_1\\
 \vdots \\ u_{n-1} \\ u_n\end{pmatrix}. \end{eqnarray*}
Note that  the differential operator $P^-$ is obtained by $P$ and the reflection $\varsigma$.  Clearly, $\mathcal{P}$ is still a linear elliptic  differential operator on $\mathbb{R}^n$.

 Let $\mathbf{K}(t,x,y)$ be the fundamental solution of the parabolic  system $\frac{\partial \mathbf{u}}{\partial t} + \mathcal{P}\mathbf{u}=0$ in $\mathbb{R}^{n}$. That is, for any $x,y\in \mathbb{R}^{n}$,
\begin{eqnarray}\label{2020.10.28-10}
   \left\{\!\begin{array}{ll} \frac{\partial \mathbf{K}(t,x,y)}{\partial t} + \mathcal{P}_{y}\mathbf{K}(t,x,y)=0,\\
    \mathbf{K}(0,x,y)=\boldsymbol{\delta}(x-y).\end{array}\right.
\end{eqnarray}
Here we write the operator $\mathcal{P}$ as $\mathcal{P}_{y}$ to denote the corresponding partial derivatives acting on variable $y$.
 Clearly, the coefficients occurring in $\mathcal{P}$ jump as $x$ crosses the $\partial \mathbb{R}^n_+$, and $\frac{\partial \mathbf{u}}{\partial t} +\mathcal{P} \mathbf{u}=0$ still has a nice fundamental solution $\mathbf{K}$ of class $C^1((0,+\infty)\times \mathbb{R}^n \times  \mathbb{R}^n) \cap C^\infty((0,+\infty)\times (\mathbb{R}^n \setminus \partial \mathbb{R}^n_+) \times (\mathbb{R}^n\setminus \partial \mathbb{R}^n_+))$, approximable even on $\partial \mathbb{R}^n_+$  by Levi's sum (see \cite{Liu-21}, or another proof below). Now, let us restrict $x,y\in \mathbb{R}^{n}_{+}$.

  Write  $\mathbf{K} (t,x,y)=\big(K^{jk} (t,x,y)\big)_{n\times n}$,  $ \;\mathbf{K}_j(t, x, y)= (K^{j1}(t,x,y)$, $\cdots$, $K^{jn} (t,x,y)$,  $\;(j=1,\cdots, n)$,   and denote   \begin{eqnarray}\label{23.12.13-2} \mathcal{F} \mathbf{K} (t,x,y) =\big( \mathcal{F} \mathbf{K}_1 (t,x,y), \cdots,  \mathcal{F} \mathbf{K}_n(t,x,y)\big)  \end{eqnarray}  
 for all $t>0,  x\in \mathbb{R}^{n}_{+}, y\in\partial \mathbb{R}^{n}_{+}$.  More precisely,  for each $j=1,\cdots, n$, the $k$-th component of $\mathcal{F} \mathbf{K}_j$ is 
 \begin{eqnarray} \label{23.12.13-6}  (\mathcal{F} \mathbf{K}_j)^k = \mu \sum_{l=1}^n (  \frac{\partial  K^{jk}}{\partial y_l} + \frac{\partial  K^{jl}}{\partial y_k})\nu_l
  +\lambda \nu_k \sum_{l=1}^n \frac{\partial  K^{jl}}{\partial y_l}.\end{eqnarray}
 It is easy to see that $\frac{\partial \big(\mathbf{K}(t,x,y)+\mathbf{K}(t,x,\overset{*}{y})\big)}{\partial \nu_y}\big|_{\partial \mathbb{R}^{n}_{+}}=
 \frac{\partial \big(\mathbf{K}(t,x,y)+\mathbf{K}(t,x,\overset{*}{y})\big)}{\partial y_n}\big|_{\partial \mathbb{R}^{n}_{+}}=
 0$
for  all $t>0$,  $x\in \mathbb{R}^{n}_{+}$ and $y\in\partial \mathbb{R}^{n}_{+}$.  Changing all  terms $\frac{\partial K^{jk}(t,x,y)}{\partial y_n}\big|_{\partial \mathbb{R}^{n}_{+}}$ into $0$ in the expression  of $2\mathcal{F}\mathbf{K}(t,x,y)$ (i.e., replacing  all terms $\frac{\partial K^{jk}(t,x,y)}{\partial y_n}\big|_{\partial \mathbb{R}^{n}_{+}}$ by $0$ in the expression of   $2\mathcal{F}\mathbf{K}(t,x,y)$),  we obtain  a matrix-valued function $\boldsymbol{\Upsilon}(t,x,y)$  for $t>0$, $x\in \mathbb{R}^{n}_{+}$ and $y\in \partial \mathbb{R}^{n}_{+}$.  This implies that  $\boldsymbol{\Upsilon}(t,x,y)$ only contains the (boundary) tangent derivatives of $K^{jk}(t, x, y)$ with respect to $y\in \partial \mathbb{R}^{n}_{+}$ (without normal derivative of $K^{jk}(t,x,y)$) in the expression). That is,  in local boundary normal coordinates (the inner normal $\boldsymbol{\nu}$ of $\partial \mathbb{R}^{n}_{+}$ is in the direction of  $x_n$-axis), \begin{eqnarray*}  \boldsymbol{\Upsilon} (t,x,y)= \Big( \boldsymbol{\Upsilon}_1 (t,x,y), \cdots,   \boldsymbol{\Upsilon}_n (t,x,y)\Big), \end{eqnarray*}
   \begin{eqnarray*}\!\!\!\!\!\!\! \!&\!\!\!&\!\!\! \boldsymbol{\Upsilon}_j (t,x,y) =
  2\mu \!\begin{small}   \begin{pmatrix} \! 2\frac{\partial K^{j1}}{\partial x_1} &\! \!\!\cdots\!  \!  & \!\!\!\frac{\partial K^{j1}}{\partial x_{n\!-\!1}} \!\!+\!\!\frac{\partial K^{\!j,n\!-\!1}}{\partial x_1}
            &   \frac{\partial K^{jn}}{\partial x_1}
      \\
\!\!\! \!\cdots \!& \!\cdots \! \!&\!\!\cdots  & \cdots \\
 \!  \! \frac{\partial K^{\!j,n\!-\!1}}{\partial x_{1}}\!\!+\! \!\frac{\partial K^{j1}} {\partial x_{n\!-\!1}}&\! \cdots    &\!\!\!2 \frac{\partial K^{j,n\!-\!1}}{\partial x_{n\!-\!1}}  &  \! \frac{\partial K^{jn}}{\partial x_{n\!-\!1}}
   \\
  \! \frac{\partial K^{jn}}{\partial x_1} &\!\!\! \cdots  \! \! & 
\frac{\partial K^{jn}}{\partial x_{n\!-\!1}}  & 
  0
 \end{pmatrix}\end{small} \!\!\begin{pmatrix} \nu_1 \\ \vdots \\ \nu_n \end{pmatrix}\\  [3mm]
\!\!\!\!\!\!&\!\!\!&\qquad \qquad \quad \; \,+ \, 2\lambda\bigg( \frac{\partial K^{j1}}{\partial x_1}   + \cdots + \frac{\partial K^{j,n-1}}{\partial x_{n-1}}  \bigg) \begin{pmatrix} \nu_1 \\ \vdots \\ \nu_n \end{pmatrix},\;\;\;  \; j=1,\cdots, n.\end{eqnarray*}
 It is easy to see that  $\boldsymbol{\Upsilon} (t,x,y)$  is a continuous (matrix-valued) function for all $t> 0$, $x\in  \mathbb{R}^{n}_{+}$ and $y\in \partial \mathbb{R}^{n}_{+}$. Further, for any  fixed $x\in  \mathbb{R}^{n}_{+}$ and any $y\in\partial \mathbb{R}^{n}_{+}$,  since $x\ne y$ we see that 
 \begin{eqnarray} \label{23.12.18-1}  \lim\limits_{t\to 0^+} \mathbf{K} (t,x,y)=0,\end{eqnarray} 
 which implies \begin{eqnarray} \label{23.12.18-2} \lim\limits_{t \to 0^{+}} \frac{\partial \mathbf{K}(t,x,y)}{\partial y_l} \big|_{\partial \mathbb{R}^{n}_{+}} =0\;\;\;\mbox{for}\;\;\,   x\in \mathbb{R}^{n}_{+}, \;\;\, y\in\partial \mathbb{R}^{n}_{+},  \;\;\, 1\le l < n. \end{eqnarray} 
     Therefore   \begin{eqnarray} \label{23.12.18-4} \lim\limits_{t\to 0^{+}} \boldsymbol{\Upsilon}(t,x,y)=0  \;\;\,\mbox{for any } x\in  \mathbb{R}^{n}_{+}, 
   \; y\in\partial \mathbb{R}^{n}_{+},\end{eqnarray}
    where $\nabla_l \mathbf{K}= (\nabla_l \mathbf{K}_1, \cdots, \nabla_l \mathbf{K}_n)$ and $\nabla_l {K}^{jk}:= \frac{\partial 
       K^{jk} (t,x,y)}{\partial y_l}  $,  $\,(1\le j,k,l\le  n)$.
     Let $\mathbf{H}(t,x,y)$ be the solution of  
    \begin{eqnarray*} \left\{ \begin{array}{ll}  \frac{\partial \mathbf{u} (t,x,y)}{\partial t} = P_g \mathbf{u} (t,x,y) \; \;\;\mbox{for}\;\,  t>0, \; x,y\in \mathbb{R}^{n}_{+},\\
2  \mu  \big(\mbox{Def}\, \mathbf{u}(t,x,y)\big)^\#  \,\boldsymbol{\nu}  +\lambda \big(\mbox{div}\, (\mathbf{u}(t,x,y))\big)\, \boldsymbol{\nu} =\boldsymbol{\Upsilon}(t,x,y)\,  \;\;\mbox{for} \,\;  t>0, \, x\in  \mathbb{R}^{n}_{+}, \, y\in \partial \mathbb{R}^{n}_{+},\\
  \mathbf{u}(0, x,y) = \mathbf{0} \,\;\;\mbox{for} \,\;  x, y \in  \mathbb{R}^{n}_{+}.\end{array} \right. \end{eqnarray*}
  From (\ref{23.12.18-4}), we get that the above parabolic system satisfy the compatibility condition.   
 Thus,  the  matrix-valued solution $\mathbf{H}(t,x,y)$, is  smooth in $(0, \infty) \times   \mathbb{R}^{n}_{+} \times  \mathbb{R}^{n}_{+}$ and  continous on  $[0, \infty) \times   {\mathbb{ R}}^{n}_{+} \times   {\mathbb{R}}^{n}_{+}$.  Then there exists a constant $C_\Omega>0$ such that  for arbitrary bounded domain $\Omega\subset \mathbb{R}^{n}_{+}$, 
  \begin{eqnarray*}  | \mathbf{H} (t,x,y) |\le C_\Omega \,\, \;\;\mbox{for all} \;\;   0\le t\le 1, \; x,\,y\in   \mathbb{R}^{n}_{+}, \;\;\,\end{eqnarray*}
   and hence, for dimensions $n\ge 2$ \begin{eqnarray}\label{23.12.9-1} \int_{\Omega} \mbox{Tr}\; \mathbf{H} (t,x,x)\, dx =nC_\Omega\,\mbox{vol}(\Omega) =o(t^{-\frac{n-1}{2}})\,  \;\;\mbox{as}\;\; t\to 0^+.\end{eqnarray}
(Actually, we have  $\,\lim_{t \to 0^{+}} \!\int_{\Omega} \mbox{Tr}\; \mathbf{H} (t,x,x)\, dx=0 $). 
   It  can easily be verified  that 
\begin{eqnarray*} 
 && \mathbf{K}^{-} (t,x,y)= \mathbf{K}(t,x,y)- \mathbf{K}(t,x,\overset{*}{y}), \\
&& \mathbf{K}^{+} (t,x,y)= \mathbf{K}(t,x,y)+ \mathbf{K}(t,x,\overset{*}{y})-  \mathbf{H}(t,x,y) ,
\end{eqnarray*}   
 are the Green functions of \begin{eqnarray*}\left\{ \begin{array}{ll}\frac{\partial \mathbf{u}}{\partial t} +{P}_g\mathbf{u}=0\;\;&\mbox{in}\;\, (0,+\infty) \times\mathbb{R}^{n}_{+},\\
\mathbf{u}=\boldsymbol{\phi} \;\; &\mbox{on}\;\; \{0\}\times \mathbb{R}^{n}_{+}\end{array}\right.\end{eqnarray*} with zero Dirichlet and zero traction (i.e., free) boundary conditions, respectively,
where $y=(y',y_n)$, $y_n\ge 0$, and  $\overset{*}{y}:=\varsigma(y',y_n)=(y^{\prime},-y_n)$.  In other words,
 \begin{eqnarray*}\left\{    \begin{array}{ll} \frac{\partial \mathbf{K}^-(t,x,y)}{\partial t} + {P}_g\mathbf{K}^-(t,x,y)=0,\;\;\; t>0,\, x,\, y\in \mathbb{R}^{n}_{+},\\
 \mathbf{K}^- (t, x,y)=0, \;\;\;\;  t>0,\; x\in  \mathbb{R}^{n}_{+}, \,\;  y\in\partial \mathbb{R}^{n}_{+},\\
 \mathbf{K}^-(0, x,y)=\boldsymbol{\delta}(x-y), \;\;\;  x,y\in \mathbb{R}^{n}_{+}\end{array} \right. \end{eqnarray*}
and
 \begin{eqnarray*}\left\{  \begin{array}{ll}   \frac{\partial \mathbf{K}^+(t,x,y)}{\partial t} + {P}_g\mathbf{K}^+(t,x,y)=0, \,\;\;\; t>0, \;x, \,y\in \mathbb{R}^{n}_{+},\\
 \mathcal {F}(\mathbf{K}^+ (t, x,y))=0, \;\; \,\,t>0,\;  x\in \mathbb{R}^{n}_{+}, \;\;  y\in\partial \mathbb{R}^{n}_{+},\\
 \mathbf{K}^+(0, x,y)=\boldsymbol{\delta}(x-y),\;\;\;  x,y\in \mathbb{R}^{n}_{+},\end{array} \right. \end{eqnarray*}
 where $ \mathcal {F}(\mathbf{K}^+ (t, x,y))$ is the traction of  $\mathbf{K}^+ (t, x,y) $ on $\partial \mathbb{R}^{n}_{+}$.
   By combining  the fact  $\frac{\partial (\mathbf{K} (t,x,y)+ \mathbf{K}(t,x,\overset{*}{y}))}{\partial y_n}\big|_{\partial \mathbb{R}^{n}_{+}} =0$  and $\mathcal{F} \big(\mathbf{K} (t,x,y)+\mathbf{K}(t, x, \overset{*}{y})\big)=\boldsymbol{\Upsilon} (t,x,y) =\mathcal{F} \mathbf{H} (t,x,y)$ for $t>0$,  $x\in \mathbb{R}^{n}_{+}$ and $y\in \partial \mathbb{R}^{n}_{+}$,  we  get $\mathcal{F} \mathbf{K}^+ (t,x,y) =0$ for all $t>0$, $x\in \mathbb{R}^{n}_{+}$ and $y\in \partial \mathbb{R}^{n}_{+}$. 
Now, we show that $\mathbf{K}^{\mp} (t,x,y)$ satisfy the corresponding boundary value problems.  In fact, for any $t>0$, $x,y \in \mathbb{R}^n_+$, we have
  $P_y\mathbf{K}(t,x, y)=\mathcal{P}_y\mathbf{K}(t,x, y)$,  so that
  \begin{eqnarray}\label{2022.10.29-8}\left\{\!\begin{array}{ll}
  \Big( \frac{\partial}{\partial t}+{P}_y\Big) \mathbf{K} (t,x, y)=
  \Big( \frac{\partial}{\partial t}+\mathcal{P}_y\Big) \mathbf{K} (t,x, y)=0,\, \;\,t>0, \,\,x,y \in \mathbb{R}^n_+,\\
 \mathbf{K} (0,x, y)= \boldsymbol{\delta}(x-y),\;\,\,x,y \in \mathbb{R}^n_+\end{array}\right.\end{eqnarray}
  by (\ref{2020.10.28-10}).
Noting that the  Jacobian matrix of the reflection $\varsigma$ is \begin{eqnarray*} \begin{pmatrix} 1 & \cdots & 0 & 0\\
\vdots & \ddots & \vdots & \vdots\\
0& \cdots & 1 & 0\\
0 & \cdots & 0 & -1\end{pmatrix},\end{eqnarray*} it follows from chain rule that for any fixed $t>0$ and $x\in \mathbb{R}^n_+$, and any $y=(y',y_n)\in \mathbb{R}^n_+$,
\begin{align*} &\left[{P}_y ( \mathbf{K}(t,x,\overset{*}{y}))\right]\bigg|_{\text{\normalsize evaluated at the point $y$}}\\
&=\left[{P}_y ( \mathbf{K}(t,x,\varsigma(y',y_n)))\right]\bigg|_{\text{\normalsize evaluated at the point $(y',y_n)$}}
 \\
&= \left[{P}_y (\mathbf{K}(t, x, (y',-y_n))\big)\right]\bigg|_{\text{\normalsize evaluated at the point $(y',y_n)$}}\\
& \! =\!
 \begin{small}{\left\{\!\left[\!-\mu\! \begin{pmatrix} \!\Delta_y\!\! & \!\cdots \!& 0 \!&\! 0\\
\vdots\!\! &\! \ddots \!& \vdots \!& \!\vdots\\
0\!\!&\! \cdots \!& \Delta_y \!& \!0\\
0\! \!&\! \cdots \!& 0 \!& \!\Delta_y\!\end{pmatrix}\!-\!(\lambda\!+\!\mu)\! \begin{pmatrix}\!\! \!\partial^2_{y_1 y_1}\! \!\!&\!\! \!\cdots \!\!\!&\! \! \partial^2_{y_1y_{n-1}} \!\!&\! \!\partial^2_{y_1y_{n}}\!\! \!\\
\!\!\!\vdots\!\!\!&\!\! \!\ddots \!&\!\vdots\!\!\!& \!\!\vdots\!\!\\
\!\!\partial^2_{y_{n-1} y_1}\! \!\!&\!\!\! \cdots \!&\! \partial^2_{y_{n-1}y_{n-1}}\! \!\!&\! \!\partial^2_{y_{n\!-\!1}y_{n}}\!\!\\
\!\!\partial^2_{y_n y_1} \!\!\!&\!\!\! \cdots \!&\! \partial^2_{y_ny_{n-1}}\! \!\!&\! \!\partial^2_{y_n y_{n}} \!\! \end{pmatrix}\!\right]\! \! \mathbf{K} (t, x, (y'\!,-y_n))\!\right\}\!\Bigg|_{\text{\normalsize evaluated at $\!(y'\!,y_n)$}}}\end{small}\\
\!& \! =
 \begin{small}{\left\{\!\left[\!-\!\mu\! \begin{pmatrix} \!\Delta_y\! &\! \cdots\! &\! 0\! &\! 0\\
\vdots \!&\! \ddots\! & \!\vdots\! & \!\vdots\\
0\!&\! \cdots\! &\! \Delta_y\! &\! 0\\
0 \!& \!\cdots\! &\! 0 \!&\! \Delta_y\end{pmatrix}\!-\!(\lambda\!+\!\mu)\! \begin{pmatrix} \!\partial^2_{y_1 y_1} \!& \!\cdots\! &\! \partial^2_{y_1y_{n-1}} \!& \! -\partial^2_{y_1y_{n}} \\
\!\cdots\!& \!\cdots\! &\!\cdots\!& \!\cdots\\
\partial^2_{y_{n-1} y_1} \!& \!\cdots \!&\! \partial^2_{y_{n-1}y_{n-1}} \!& \!-\partial^2_{y_{n\!-\!1}y_{n}}\\
-\partial^2_{y_n y_1} \!&\! \cdots \!&\! -\partial^2_{y_ny_{n-1}} \!& \!\partial^2_{y_n y_{n}}  \end{pmatrix}\!\right]\!  \mathbf{K}\!\right\}\Bigg|_{\text{\normalsize evaluated  at $(y',-y_n)$}}}\end{small} \\
&= P^-_{\overset{*}{y}} ( \mathbf{K}(t, x, \overset{*}{y}))\Big|_{\text{\normalsize evaluated at the point $\overset{*}{y}=(y',-y_n)$}}.  \end{align*}
 That is, the action of $P_y$ to $\mathbf{K}(t,x,\overset{*}{y})$ at the point ${y}=(y',y_n)$ is just the action of $P^-_{\overset{*}{y}} $
 to $\mathbf{K}(t,x,\overset{*}{y})$ at the point $\overset{*}{y}=(y',-y_n)$. Because of  $\varsigma(y',y_n)=(y',-y_n)\in \mathbb{R}^n_-$, we see $$P^-_{\overset{*}{y}}( \mathbf{K}(t,x,\overset{*}{y}))\big|_{\text{\normalsize evaluated at the point $\overset{*}{y}=(y',-y_n)$}}=\mathcal{P}_{\overset{*}{y}}(\mathbf{K} (t,x,\overset{*}{y}))\big|_{\text{\normalsize evaluated at the point $\overset{*}{y}=(y',-y_n)$}}.$$
For any $t>0$, $x\in \mathbb{R}^n_+$ and $(y',-y_n)\in \mathbb{R}^n_-$, we have  $$(\frac{\partial}{\partial t}+ {\mathcal{P}_y}) (\mathbf{K} (t, x,  (y',-y_n)))=0.$$
In addition, $\mathbf{K} (t, x,  (y',-y_n))= \mathbf{K} (t, x,  \varsigma(y))$ for any $t>0$, $x,y\in \mathbb{R}^n_+$.
By virtue of  $x\ne (y',-y_n)$, this leads to
 \begin{eqnarray*} \label{2020.10.29-1} \left\{ \begin{array}{ll} (\frac{\partial}{\partial t}+ {{P}_y}) \big(\mathbf{K} (t, x,  (y',-y_n))\big)=0\;\; \mbox{for any}\;\,t>0, x\in \mathbb{R}^n_+ \;\,\mbox{and}\;\,(y',-y_n)\in \mathbb{R}^n_-,\\
 \mathbf{k}(0,x,(y',-y_n))=0 \;\,\mbox{for any}\;\, x\in \mathbb{R}^n_+ \;\,\mbox{and}\;\,(y',-y_n)\in \mathbb{R}^n_-,\end{array}\right.\end{eqnarray*}
 i.e.,
  \begin{eqnarray} \label{2020.10.29-2}\left\{\begin{array}{ll}  (\frac{\partial}{\partial t}+ {{P}_y}) \mathbf{K} (t, x,  \overset{*}{y})=0\;\; \mbox{for any}\;\,t>0, x\in \mathbb{R}^n_+ \;\,\mbox{and}\;\,\overset{*}{y}\in \mathbb{R}^n_-\\
  \mathbf{k}(0,x,\overset{*}{y})=0 \;\,\mbox{for any}\;\, x\in \mathbb{R}^n_+ \;\,\mbox{and}\;\,\overset{*}{y}\in \mathbb{R}^n_-.\end{array}\right.\end{eqnarray}
 Combining (\ref{2022.10.29-8}) and (\ref{2020.10.29-1}), we obtain that
 \begin{eqnarray} \label{2020.10.29-3} \left\{\! \begin{array}{ll}(\frac{\partial}{\partial t}+ {{P}_y}) \Big(\mathbf{K} (t, x, {y})-\mathbf{K} (t, x,  \overset{*}{y})\Big)=0\;\;\,\mbox{for any}\;\,t>0,\;\, x,y\in \mathbb{R}^n_+\\
 \mathbf{K} (0, x, {y})-\mathbf{K} (0, x,  \overset{*}{y})=\boldsymbol{\delta} (x-y)\;\;\, \mbox{for any}\;\, x,y\in \mathbb{R}^n_+\end{array}\right.
 .\end{eqnarray}
  $\mathbf{K}(t,x, y)$ is $C^1$-smooth with respect to $y$ in $\mathbb{R}^n$ for any fixed $t>0$ and $x\in \mathbb{R}^n_+$, so does it on the hyperplane $\partial \mathbb{R}^n_+$. Therefore, we get that  $\mathbf{K}^-(t,x,y)$ (respectively $\mathbf{K}^+(t,x,y)$) is the Green function in the upper-half space with the Dirichlet (respectively, Neumann) boundary condition on $\partial \mathbb{R}^n_+$.

\vskip 0.12 true cm

  To show $C^1$-regularity of the fundamental solution $\mathbf{K}(t,x,y)$, it suffices to prove  $C^{1,1}_{loc}$-regularity for the solution $\mathbf{u}(t,x)$ of the  parabolic system $\frac{\partial \mathbf{u}(t,x)}{\partial t} + \mathcal{P}\mathbf{u}(t,x)=0$  in $(0,+\infty)\times \mathcal{M}$. This immediately follows from Dong's celebrated result \cite{Do-12} of local $C^{1,1}$-regularity for the parabolic equation 
  \begin{eqnarray*} \left\{ \begin{array}{ll} \frac{\partial{u}}{\partial t}+ L {u}  =f \;\;\;\mbox{in}\;\; (0, +\infty)\times U,\\
  u=\phi \;\;\;\mbox{on}\;\;\, \{0\}\times \partial U,\end{array} \right.\end{eqnarray*}
 where $L$ is an elliptic differential operator whose coefficients and data  are irregular in one spatial direction (see Remark 5 of p. 141 in \cite{Do-12}).   In fact,  the coefficients of $\mathcal{P}$ are smooth in $\mathcal{M}\setminus (\partial \Omega)$, and the coefficients of $\mathcal{P}$  is irregular only in the normal direction  (but the coefficients of $\mathcal{P}$ are all smooth in other  direction) on $\partial \Omega$.  
  \vskip 0.25 true cm
Generally, let ${L}={L}(\partial_{x_1},\cdots, \partial_{x_n})$ be a  linear differential system with respect to $\partial_{x_1}$, $\cdots$,  $\partial_{x_n}$, which is defined on $\mathbb{R}^n_+$. By reflecting (or doubling) we can extend the $L$ to the whole Euclidean space $\mathbb{R}^n$ as follows:  on the upper half-space $\mathbb{R}^n_+$, the differential system is still  ${L}(\partial_{x_1}, \cdots, \partial_{x_n})$; but on the lower half-space $\mathbb{R}^n_-$, the  differential system is modified as ${L}(\partial_{x_1},\cdots, \partial_{x_{n-1}}, -\partial_{x_n})$ (i.e., a reflection differential system). We define $\mathcal{L}$ as
\begin{eqnarray} \mathcal{L}:= \left\{\begin{array}{ll}  {L}(\partial_{x_1}, \cdots, \partial_{x_n})\;\;\;\;\,\;\;\; \mbox{on the upper half-space}\;\, \mathbb{R}^n_+,\\
 {L}(\partial_{x_1}, \cdots, -\partial_{x_n})\;\;\;\;\;\mbox{on the lower half-space}\;\, \mathbb{R}^n_-.\end{array}\right.\end{eqnarray}
 Such an extension bases on the reflection $\varsigma$ with respect to $x'$-hyperplane, which has changed the representation  of differential system on the lower half-space from the given differential system $L$ (defined on the upper half-space). By a completely similar method, we can also get its Green's function with corresponding Dirichlet (or Neumann) boundary condition in the upper half-space $\mathbb{R}^n_+$. Furthermore, such a basic ``reflection method'' for obtaining  Green's function of a parabolic (elastic) system in $\mathbb{R}^{n}_{+}$ with zero Dirichlet (or Neumann) boundary condition can be generalized to a smooth Riemannian manifold with smooth boundary, just as done in \cite{Liu-21} for the Lam\'{e} operator (even for more general linear elliptic system).
 This technique is also called the ``method of image'' in \cite{Liu-21}. Let us introduce it again (the details can been seen in \cite{Liu-21}).  Let $\mathcal{M}=\Omega \cup (\partial \Omega)\cup \Omega^*$ be the (closed) double of $\Omega$.
 Note that  $\mathcal{M}$ is a closed manifold since we can take the boundary normal coordinates on a neighborhood of $\partial \Omega$, and used the metric (see \cite{Liu-21.1}, \cite{LiuTan-22} or p.\,10169, p.\,10183 and p.\,10187 of \cite{Liu-21} )
\begin{eqnarray} \label{2021.2.6-3}  g_{jk} (\overset{*}{x})\!\!\!&\!=\!&\!\!\!- g_{jk} (x) \quad \, \mbox{for}\;\;
  j<k=n \;\;\mbox{or}\;\; k<j=n,\\ \!\!\!&\!=\!&\!\!\! g_{jk} (x) \;\;\mbox{for}\;\; j,k<n \;\;\mbox{or}\;\; j=k=n,\nonumber\\
  \label{2021.2.6-4}  g_{jk}(x)\!\!\!&\!=\!&\!\!\! 0 \;\; \mbox{for}\;\; j<k=n \;\;\mbox{or}\;\; k<j=n \;\;\mbox{on}\;\; \partial \Omega,\\
   \label{2021.2.6-5}  \sqrt{|g|/g_{nn}} \; dx_1\cdots dx_{n-1}\!\!\!&\!=\!&\!\!\! \mbox{the element of (Riemannian) surface area on} \,\, \partial \Omega,\end{eqnarray} where $x_n(\overset{*}{x})= -x_n (x)$.
    This follows from the following fact: recall that the Riemannian metric $(g_{ij})$ is given in the local coordinates $x_1, \cdots, x_n$, i.e., $g_{ij}(x_1,\cdots, x_n)$. In terms of new coordinates $z_1,\cdots, z_n$, with  $x_i=x_i(z_1,\cdots, z_n), \,\, i=1,\cdots, n,$  the same metric is given by the functions $\tilde{g}_{ij} =\tilde{g}_{ij} (z_1, \cdots, z_n)$, where
\begin{eqnarray} \label{2022.10-2} \tilde{g}_{ij}= \frac{\partial x_k}{\partial z_i} g_{kl} \frac{\partial x_l}{\partial z_j}.\end{eqnarray}
  If $S$ is a local  (opposite orientation) coordinate change in a neighborhood intersecting with $x'$-hyperplane
  \begin{eqnarray}\label{2022.10.28-1} \left\{ \begin{array}{ll} x_1= z_1, \\
                   \cdots \cdots\\
                   x_{n-1}=z_{n-1},\\
                   x_n=-z_n,\end{array}\right.\end{eqnarray}
                                       then its Jacobian matrix is
\begin{eqnarray*}     \begin{pmatrix} 1& \cdots &   0 & 0\\
   \vdots & \ddots & \vdots & \vdots\\
  0 & \cdots & 1 &  0\\
   0& \cdots & 0&-1\end{pmatrix}.\end{eqnarray*}
Using this and (\ref{2022.10-2}), we immediately obtain (\ref{2021.2.6-3})--(\ref{2021.2.6-4}).
\vskip 0.15 true cm

On a smooth Riemannian manifold $(\Omega,g)$ with smooth boundary $\partial \Omega$, since the Navier-Lam\'{e} operator is a linear differential operator defined on $\Omega$, it can be further written as $P_g:=P(g^{\alpha\beta}(x)$, $g^{\alpha n}(x), g^{n\beta}(x)$, $g^{nn}(x)$, $\frac{\partial}{\partial x_1}$, $\cdots$, $\frac{\partial}{\partial x_{n-1}}$, $\frac{\partial }{\partial x_n})$.
We  define \begin{eqnarray} \label{2022.10.18-2} \mathcal{P}=\left\{\begin{array}{ll} \! P \;\;\; \;\;\;\,\mbox{on} \;\, \Omega,\\
 \!P^* \;\;\;\; \;\mbox{on} \;\, \Omega^*, \end{array} \right.\end{eqnarray}
where \begin{eqnarray} \label{2022.10.6-8} P^*:=P\Big({g}^{\alpha\beta}(\overset {*}{x}), - {g}^{\alpha n}(\overset{*}{x}),- {g}^{n\beta}(\overset{*}{x}), {g}^{nn}(\overset{*}{x}), \frac{\partial }{\partial x_1}, \cdots, \frac{\partial }{\partial x_{n-1}},- \frac{\partial }{\partial x_n}\Big), \end{eqnarray}
and $\overset{*}{x}=(x',-x_n)\in \Omega^*$.
  Clearly, the differential operator $P^*$ is obtained by $P_g$ and the reflection $\varsigma$, that is,  $P^*$ is got if we replace $g^{\alpha\beta}(x)$, $g^{\alpha n}(x)$, $g^{n\beta} (x)$, $g^{nn}(x)$, $\frac{\partial }{\partial x_n}$ by ${g}^{\alpha\beta}(\overset{*}{x})$, $-{g}^{\alpha n}(\overset{*}{x})$, $-{g}^{n\beta} (\overset{*}{x})$, ${g}^{nn}(\overset{*}{x})$, $-\frac{\partial }{\partial x_n}$ in $P\Big( g^{\alpha\beta} (x)$, $g^{\alpha n}(x)$, $g^{n\beta}(x)$, $g^{nn}(x)$, $\frac{\partial }{\partial x_1}$, $\cdots$, $\frac{\partial }{\partial x_{n-1}}$, $\frac{\partial }{\partial x_n}\Big)$, respectively.  Note that  $g^{\alpha\beta} (\overset{*}{x})=
g^{\alpha\beta} (x)$, $-g^{\alpha n}(\overset{*}{x})= g^{\alpha n}(x)$, $-g^{n\beta}(\overset{*}{x})=g^{n\beta}(x)$ and $g^{nn}(\overset{*}{x})=g^{nn} (x)$.
In view of the metric matrices $g$ and $g^*$ have the same order principal minor determinants, we see that $\mathcal{P}$ is still a linear elliptic differential operator on $\mathcal{M}$.

\vskip 0.15 true cm

    For given $\epsilon >0$, denote by $U_{\epsilon}(\partial \Omega)=\{z\in \mathcal{M}\big|\mbox{dist}(z, \partial \Omega)<\epsilon \}$ the $\epsilon$-neighborhood of $\partial \Omega$ in $\mathcal{M}$. Clearly, for any $W\subset U_\epsilon (\partial \Omega)$ with local coordinates $x=(x',x_n)$ such that $\epsilon >x_n>0$ in $W\cap \Omega$;  $x_n=0$ on $W\cap \partial\Omega$, and $|\nabla x_n|=1$ near $\partial \Omega$. It is easy to verify (see \cite{Liu-21}) that  such an extend metric $g$ is $C^0$-smooth on whole $\mathcal{M}$ and $C^\infty$-smooth in $\mathcal{M}\setminus \partial \Omega$.  Let us point out that obtaining a closed double manifold $\mathcal{M}$ from $\Omega$ is a standard and classical method in Riemannian geometry (see, Example 9.32 on p.$\,226$ in \cite{Lee-13} or \cite{MS-67}). Roughly speaking, $\Omega$ and $\Omega^*$ are $C^0$-smoothly glued into $\mathcal{M}$. Next, let $\mathcal{P}$ be the double to $\mathcal{M}$ of
 the  operator $P_g$ on $\Omega$ ($\mathcal{P}$ has new representation on $\Omega^*$ according to (\ref{2022.10.6-8}). Then the coefficients of $\mathcal{P}$ are also uniformly bounded on $\mathcal{M}$ and $C^\infty$-smooth in $\mathcal{M}\setminus \partial \Omega$ (the coefficients of  $\mathcal{P}$ are smooth up to the boundary $\partial \Omega$ from two sides of $\partial \Omega$, but the coefficients occurring in $\mathcal{P}$ jump as $x$ crosses $\partial \Omega$).
 Of course, $\mathcal{P}$ is an elliptic operator on $\mathcal{M}$.
Thus we immediately get a strongly continuous semigroup $e^{-t\mathcal{P}}$ on $\mathcal{M}$. Let $\mathbf{K}(t,x,y)$ be the parabolic kernel (i.e., a fundamental solution to parabolic  system $\frac{\partial \mathbf{u}}{\partial t} + \mathcal{P}\mathbf{u}=0$ on whole (closed) manifold $\mathcal{M}$). More precisely, \begin{align*}\left\{\!\begin{array}{ll}
           \frac{\partial \mathbf{K}(t,x,y)}{\partial t} + \mathcal{P}_y \mathbf{K}(t,x,y)=0,\quad t>0,\;\; x,y\in \mathcal{M}, \\
        \mathbf{K}(0,x,y) = \boldsymbol{\delta}(x-y),\;\; x,y\in \mathcal{M}\end{array}\right.  \end{align*}
           where we have written $\mathcal{P}$ as $\mathcal{P}_y$ to denote the corresponding covariant derivatives acting on variable $y$. Obviously, $\mathbf{K}\in C^1((0,\infty)\times \mathcal{M} \times \mathcal{M}) \cap C^\infty ((0,\infty)\times (\mathcal{M}\setminus \partial \Omega)\times (\mathcal{M}\setminus \partial \Omega))$.
           Let  $\mathbf{H}(t,x,y)$ be the solution of  
    \begin{eqnarray*} \left\{ \begin{array}{ll}  \frac{\partial \mathbf{u} (t,x,y)}{\partial t} = P_g \mathbf{u} (t,x,y) \; \;\;\mbox{for}\;\,  t>0, \; x,y\in \Omega,\\
2  \mu  \big(\mbox{Def}\, \mathbf{u}(t,x,y)\big)^\#  \,\boldsymbol{\nu}  +\lambda \big(\mbox{div}\, (\mathbf{u}(t,x,y))\big)\, \boldsymbol{\nu} =\boldsymbol{\Upsilon}(t,x,y)\,  \;\;\mbox{for} \,\;  t>0, \, x\in \Omega, \, y\in \partial \Omega,\\
  \mathbf{u}(0, x,y) = \mathbf{0} \,\;\;\mbox{for} \,\;  x, y \in \Omega, \end{array} \right. \end{eqnarray*}
where 
   \begin{eqnarray*}  \boldsymbol{\Upsilon} (t,x,y)= \Big( \boldsymbol{\Upsilon}_1 (t,x,y), \cdots,   \boldsymbol{\Upsilon}_n (t,x,y)\Big), \end{eqnarray*}
   \begin{eqnarray*}\!\!\!\!\!\!\! \!&\!\!\!&\!\!\! \boldsymbol{\Upsilon}_j (t,x,y) \\
\!\!\!\!\!\! \!\!\! &\!\!\!\!\!&\!\!\! \! =\!
  2\mu \!\begin{small}   \begin{pmatrix} \! 2\frac{\partial K^{j1}}{\partial x_1}+2\!\sum\limits_{m=1}^n \!\Gamma^{1}_{1m} K^{jm}  &\! \!\!\cdots\!  \!  & \!\!\!\frac{\partial K^{j1}}{\partial x_{n\!-\!1}} \!\!+\!\!\frac{\partial K^{\!j,n\!-\!1}}{\partial x_1} \!\!+\!\!\! \sum\limits_{m=1}^n\!\! \big( \Gamma^1_{\!n\!-\!1, m}\! \!+\!\!\Gamma^{n\!-\!1}_{\!1m} \big) K^{jm} 
            &   \frac{\partial K^{jn}}{\partial x_1}\!\!+\!\!\!\sum\limits_{m=1}^n\!\! \big( \Gamma^1_{\!nm} \!\!+\!\!\Gamma^n_{1m}\big) K^{jm}  
      \\
\!\!\! \!\cdots \!& \!\cdots \! \!&\!\!\cdots  & \cdots \\
 \!  \! \frac{\partial K^{\!j,n\!-\!1}}{\partial x_{1}}\!\!+\! \!\frac{\partial K^{j1}} {\partial x_{n\!-\!1}}\! \!+\! \!\!\sum\limits_{m=1}^n \!\! \big( \Gamma^{n\!-\!1}_{\!1m}\!\! +\!\! \Gamma^1_{\!n\!-\!1, m} \big)  K^{jm}
 \! \!&\! \cdots    &\!\!\!2 \frac{\partial K^{j,n\!-\!1}}{\partial x_{n\!-\!1}}\!\!+\!\!2\!\sum\limits_{m=1}^n \!\!\Gamma^{n\!-\!1}_{\!n\!-\!1,m} K^{jm} &  \! \frac{\partial K^{jn}}{\partial x_{n\!-\!1}}\!\!+\!\!\sum\limits_{m=1}^n\!\! \big(\Gamma^{n\!-\!1}_{\!\!nm}\! \!+\!\!\Gamma_{\!\!n\!-\!1,m}^n \!\big) K^{\!jm} \!\!
   \\
  \! \frac{\partial K^{jn}}{\partial x_1} \!+\!\!\sum\limits_{m=1}^n \!\!\big( \Gamma_{\!1m}^n\! +\!\Gamma^1_{\!nm}\big) K^{jm}   \! &\!\!\! \cdots  \! \! & 
\frac{\partial K^{jn}}{\partial x_{n\!-\!1}}\! +\!\big(\Gamma^n_{n-1,m}+\Gamma_{\!nm}^{n\!-\!1} \big)K^{jm}   & 
    \;2 \sum\limits_{m=1}^n 
 \Gamma^n_{nm}  K^{jm}\!\!
 \end{pmatrix}\end{small} \!\!\begin{pmatrix} \nu_1 \\ \vdots \\ \nu_n \end{pmatrix}\\  [3mm]
\!\!\!\!\!\!&\!\!\!&\!\!\! \;\;\; + \, 2\lambda\bigg( \frac{\partial K^{j1}}{\partial x_1} + \sum_{m=1}^n \Gamma^{1}_{1m} K^{jm}  + \cdots + \frac{\partial K^{j,n-1}}{\partial x_{n-1}} + \sum_{m=1}^n \Gamma^{n-1}_{n-1,m} K^{jm} +  \sum_{m=1}^n \Gamma^{n}_{nm} K^{jm} \bigg) \begin{pmatrix} \nu_1 \\ \vdots \\ \nu_n \end{pmatrix},\;\;\;  \; j=1,\cdots, n.\end{eqnarray*}        
       From (\ref{23.12.18-4}), we get that the above parabolic system satisfy the compatibility condition.   
 Thus,  the  matrix-valued solution $\mathbf{H}(t,x,y)$, is  smooth in $(0, \infty) \times  \Omega \times  \Omega$ and  continous on  $[0, \infty) \times \bar \Omega \times \bar \Omega$.  Then there exists a constant $C>0$ such that 
  \begin{eqnarray*}  | \mathbf{H} (t,x,y) |\le C \,\, \;\;\mbox{for all} \;\;   0\le t\le 1, \; x,\,y\in \bar \Omega, \;\;\,\end{eqnarray*}
   and hence, for dimensions $n\ge 2$,   \begin{eqnarray}\label{23.12.9-1} \int_{\Omega} \mbox{Tr}\; \mathbf{H} (t,x,x)\, dx =nC\,\mbox{vol}(\Omega) =o(t^{-\frac{n-1}{2}})\,  \;\;\mbox{as}\;\; t\to 0^+.\end{eqnarray}
(Actually, we have  $\,\lim_{t \to 0^{+}} \!\int_{\Omega} \mbox{Tr}\; \mathbf{H} (t,x,x)\, dx=0 $). 
  Therefore, \begin{eqnarray} \label{2022.7.31} \mathbf{K}^-(t,x, y) = \mathbf{K}(t,x,y) -\mathbf{K}(t,x, \overset{*}{y}), \;\; \, x,y\in \Omega \end{eqnarray}
and \begin{eqnarray}\label{2022.8.1} \mathbf{K}^+(t,x, y) = \mathbf{K}(t,x,y) +\mathbf{K}(t,x, \overset{*}{y})-\mathbf{H}(t,x,y), \;\; \, x,y\in \Omega \end{eqnarray}
are Green's functions of the parabolic (elastic) system $\frac{\partial\mathbf{u} }{\partial t} + P_g\mathbf{u}=0$ in $\Omega$ with zero Dirichlet and Neumann  boundary conditions, respectively.    
That is,
\begin{align} \label{2022.7.30-5}
    \begin{cases}
        \frac{\partial \mathbf{K}^{-}(t,x,y)}{\partial t} + P_g \mathbf{K}^{-}(t,x,y)=0,\quad t>0,\; x,y\in \Omega,\\
        \mathbf{K}^{-}(t,x,y)=0,\quad t>0,\; x\in\Omega,\ y\in\partial\Omega,\\
      \mathbf{K}^{-}(0,x,y)=  \boldsymbol{\delta}(x-y), \quad x,y\in \Omega
    \end{cases}
\end{align}
\begin{align}\label{2022.7.30-6}
    \begin{cases}
        \frac{\partial \mathbf{K}^{+}(t,x,y)}{\partial t} + P_{g}\mathbf{K}^{+}(t,x,y)=0,\quad t>0,\; x,y\in \Omega,\\
        \frac{\partial \mathbf{K}^{+}(t,x,y)}{\partial \boldsymbol{\nu}}=0,\quad t>0,\; x\in\Omega,\, y\in\partial\Omega,\\
       \mathbf{K}^{+}(0,x,y) =  \boldsymbol{\delta}(x-y), \quad x,y\in \Omega
    \end{cases}
\end{align}
where $\frac{\partial \mathbf{K}^{+}}{\partial \boldsymbol{\nu}}:=2\mu(\operatorname{Def}\,\mathbf{K}^{+})^\#\boldsymbol{\nu}+\lambda\,(\operatorname{div}\mathbf{K}^{+})\boldsymbol{\nu}$ is the Neumann boundary condition again. Note that all covariant derivatives act on variable $y$ for any fixed $x\in \Omega$. This method based on theory of parabolic system on Riemannian manifolds. By the ``method of image'' and (\ref{2022.7.31})--(\ref{2022.7.30-6}), the studied topic becomes to discuss the differences of the two parts on $\Omega$ and $\Omega^*$ for a fundamental solution $\mathbf{K}$ on $\mathcal{M}$, further becomes to study such a problem which is restricted on $\Omega$. Of course, one can get all the coefficients
$a_l^\mp$ by classical estimates (it needs a great deal of arduous estimates and and tedious calculations). But we have used the method of pseudodifferential operators in \cite{Liu-21} (which is quite simpler and more effective) as follows.

It is well-known that \begin{eqnarray} \label{1-0a-2}\int_\Omega  \big(\mbox{Tr}\,({\mathbf{K}}^\mp(t,x,x))\big) dV=\sum_{k=1}^\infty e^{-t \tau_k^\mp}.\end{eqnarray}
  Since $e^{-t\mathcal{P}}f(x) =\frac{1}{2\pi i} \int_{\mathcal{C}} e^{-t\tau } (\tau I- \mathcal{P})^{-1} f(x) \, d\tau$, we have
  $$e^{-t\mathcal{P}} f(x) = \frac{1}{(2\pi)^n} \int_{{\mathbb{R}}^n} e^{ix\cdot \xi} \Big( \frac{1}{2\pi i} \int_{\mathcal{C}} e^{-t\tau}\, \iota\big((\tau I-\mathcal{P})^{-1}\big) \, \hat{f} (\xi) \, d\tau\Big) d\xi,$$
  so that \begin{eqnarray*}&&\mathbf{K}(t,x,y)= e^{-t\mathcal{P}} \boldsymbol{\delta}(x-y) = \frac{1}{(2\pi)^n} \int_{{\mathbb{R}}^n} e^{i(x-y)\cdot \xi} \Big( \frac{1}{2\pi i} \int_{\mathcal{C}}  e^{-t\tau}\, \iota\big((\tau I-\mathcal{P})^{-1}\big) \, d\tau\Big) d\xi\\
  && \qquad \qquad  \; = \frac{1}{(2\pi)^n} \int_{{\mathbb{R}}^n} e^{i(x-y)\cdot \xi} \Big( \frac{1}{2\pi i} \int_{\mathcal{C}}  e^{-t\tau}\, \big(\sum_{l\ge 0} {\mathbf{q}}_{-2-l} (x, \xi,\tau)\big) \, d\tau\Big) d\xi,\end{eqnarray*}
  where $\mathcal{C}$ is a suitable curve in the complex plane in the positive direction around the spectrum of $\mathcal{P}$, and $\iota((\tau I-\mathcal{P})^{-1}):=\sum_{l\ge 0} {\mathbf{q}}_{-2-l} (x,\xi,\tau)$ is the full symbol of resolvent operator $(\tau I-\mathcal{P})^{-1}$.
  This implies that for any $x\in \bar \Omega$,
  \begin{eqnarray*}&&\mbox{Tr}\,( {\mathbf{K}}(t,x,x))= \frac{1}{(2\pi)^n} \int_{{\mathbb{R}}^n}   \Big( \frac{1}{2\pi i} \int_{\mathcal{C}}  e^{-t\tau}\,\sum_{l\ge 0}  \mbox{Tr}\,({\mathbf{q}}_{-2-l} (x,\xi, \tau) )\, d\tau\Big) d\xi,\\
  &&  \mbox{Tr}\,({\mathbf{K}}(t,x,\overset{*}{x}))= \frac{1}{(2\pi)^n} \int_{{\mathbb{R}}^n}  e^{i(x-\overset{*}{x})\cdot \xi}   \Big(  \frac{1}{2\pi i} \int_{\mathcal{C}}  e^{-t\tau}\,\sum_{l\ge 0} \mbox{Tr}\,({\mathbf{q}}_{-2-l} (x,\xi, \tau)) \, d\tau\Big) d\xi
  .\end{eqnarray*}
  It is easy to see that  the traces of the symbols (i.e., $\mbox{Tr} \, (\mathbf{K} (t,x,y))$ and $\mbox{Tr} \, (\mathbf{K} (t,x,\overset{*}{y}))$),
  have the same expressions either on $\Omega$ or on $\Omega^*$.

 Therefore, whether one can get all coefficients of the parabolic trace, it depends strongly on whether one can explicitly calculate the full symbol of the resolvent operator $(\tau I- \mathcal{P})^{-1}$. The second key idea (or main contribution) of \cite{Liu-21} is that such kind of full symbol can be explicitly obtained by a surprising new technique \cite{Liu-19} (see below).

\vskip 0.12 true cm

 Since the full symbol $\sigma((\tau I-P_g)^{-1})$ of $(\tau I-P_g)^{-1}$ is classical, we write the full symbol $\sigma((\tau I-P_g)^{-1})\sim \mathbf{q}_{-2}+\mathbf{q}_{-3}+\mathbf{q}_{-4} +\cdots$ with $\mathbf{q}_j$ being homogeneous of degree $j\ (j=-2,-3,\cdots)$ in the variable $(\xi, \tau^{1/2})$. Generally, for an elliptic system, it is not possible to explicitly calculate the full symbol of its resolvent operator. However, for the elastic Lam\'{e} system, we discovered that the full symbol is explicitly calculable. In \cite{Liu-21}, we have explicitly calculated the full symbol by a method of operator algebra (This is a surprising effective method stemmed from the idea of \cite{Liu-19} for studying the elastic Dirichlet-to-Neumann map). Thus by the parabolic trace method, we can get all coefficients  $a_0$, $a_1^\mp$, $a_2^\mp$, $\cdots$
in the expansion of the parabolic trace
\begin{eqnarray*} \sum_{k=1}^\infty e^{-t \tau_k^\mp} = t^{-n/2} \left[ a_0 + a_1^\mp t^{1/2} +a^\mp_3 t + a^\mp_4 t^{3/2} +\cdots\right].\end{eqnarray*}
In \cite{Liu-21}, we had obtained  $a_0$ and $a_1^\mp$ by an (algorithm) method. From the expression (3.10)  in \cite{Liu-21}, we can get $\mathbf{q}_{-3}$ and $\mathbf{q}_{-4}$, i.e.,
\begin{eqnarray} \label{2022.8.2-1} {\mathbf{q}}_{-3}(x,\xi,\tau)= -{\mathbf{a}}_2^{-1}\Big({\mathbf{a}}_1  {\mathbf{a}}_2^{-1} - i \sum\limits_{l=1}^n \frac{\partial {\mathbf{a}}_2}{\partial \xi_l}\frac{\partial {\mathbf{a}}_2^{-1}}{\partial x_l} \Big),\\
{\mathbf{q}}_{-4} = -{\mathbf{a}}_2^{-1} \bigg( \sum_{\substack{j<2\\ j+|\alpha|-k=0}} (\partial^\alpha_\xi  {\mathbf{a}}_k)(D^\alpha_x {\mathbf{q}}_{-2-j})/\alpha !\bigg). \end{eqnarray}
 Note that $\mathbf{q}_{-3}$ has not any contribution to the coefficients $a_l^\mp$ because it is an odd function in variable $\xi\in \mathbb{R}^n$ (the corresponding integral on $\mathbb{R}^n$ will vanish). Further, by inserting  $\mathbf{q}_{-4}$ into our calculation procedure for the coefficients $a_2^\mp$ (see \cite{Liu-21}), we immediately get the explicit expression of $a_2^\mp$ (i.e., (\ref{2022.7.31-6})). To this end, we need to use the expansion of the metric $g$:
 we can write
          \begin{eqnarray} \label{5,,1}   g_{jk} (x) \!\! &=&\!\!  \delta_{jk} +\sum_{l=1}^{n} \frac{\partial g_{jk}}{\partial x_{l}}(0) \, x_{l}+
          \frac{1}{2} \sum_{l,m=1}^{n} \frac {\partial^2 g_{jk} }{\partial x_l\partial x_m} (0)\, x_l x_m +O(|x|^3)  \; \; \mbox{near}\;\, 0. \nonumber\end{eqnarray}
In geodesic
normal coordinate system centered at $x=0$, as Riemann showed, one has
\begin{eqnarray*} \frac{\partial^2g_{lm}}{\partial x_j\partial x_j}=-\frac{1}{3}R_{ljmk}-\frac{1}{3} R_{ lkmj}.\end{eqnarray*}
Substituting this into $\mathbf{q}_{-4}$ in \cite{Liu-21}, and by calculating $$\int_{ W\cap \Omega} \left\{ \frac{1}{(2\pi)^n} \int_{\mathbb{R}^n} \Big( \frac{1}{2\pi i} \int_{\mathcal{C}} e^{-t\tau} \, \mbox{Tr}\,\big( \mathbf{q}_{-4} (x, \xi, \tau)\big) d\tau\Big) d\xi\right\} dV,$$ where $W\subset \Omega\setminus U_\epsilon (\partial \Omega)$,
we obtain the term that will product the factor $\frac{t}{3}\int_\Omega R\, dV$.
Note that  under the boundary normal coordinates, $\frac{1}{2}\frac{\partial g_{jk}}{\partial x_n}= \delta_{jk} \tilde{\kappa}_j$ on $\partial \Omega$. By calculating $$\int_{\tilde{W}\cap \Omega} \left\{ \frac{1}{(2\pi)^n} \int_{\mathbb{R}^n} \Big( \frac{1}{2\pi i} \int_{\mathcal{C}} e^{-t\tau} \, \mbox{Tr}\,\big( \mathbf{q}_{-4} (x, \xi, \tau)\big) d\tau\Big) d\xi\right\} dV,$$ we get  another term that  will product the factor $-\frac{t}{6}\int_{\tilde{W}\cap \partial \Omega} H\, ds$, where $\tilde{W}\subset U_\epsilon (\partial \Omega)$.
In addition, it is easy to verify that
\begin{eqnarray*} && \int_{{W}\cap \Omega} \left\{ \frac{1}{(2\pi)^n} \int_{\mathbb{R}^n} e^{i(x-\overset{*}{x})\cdot \xi} \Big( \frac{1}{2\pi i} \int_{\mathcal{C}} e^{-t\tau} \, \mbox{Tr}\,\big( \mathbf{q}_{-4} (x, \xi, \tau)\big) d\tau\Big) d\xi\right\} dV= O(t^{2-\frac{n}{2}}) \;\;\; \mbox{as}\;\; t\to 0^+,\\
&& \int_{\tilde{W}\cap \Omega} \left\{ \frac{1}{(2\pi)^n} \int_{\mathbb{R}^n} e^{i(x-\overset{*}{x})\cdot \xi} \Big( \frac{1}{2\pi i} \int_{\mathcal{C}} e^{-t\tau} \, \mbox{Tr}\,\big( \mathbf{q}_{-4} (x, \xi, \tau)\big) d\tau\Big) d\xi\right\} dV= O(t^{2-\frac{n}{2}}) \;\;\; \mbox{as}\;\;\; t\to 0^+.\end{eqnarray*}
 Therefore, Theorem 1.2 is proved.

Based on classical theory of strongly continuous semigroups and pseudodifferential operators, we had given a rigorous proof to Theorem 1.1 in  \cite{Liu-21}. Here we also give an outline of proof of Theorem 1.2. Clearly, the conclusion and proof of \cite{Liu-21} are correct.

\vskip 0.86 true cm

\section{Several remarks}

 \vskip 0.26 true cm

\begin{remark}
    Besides the elastic Lam\'{e} system, by applying our new technique we can also explicitly calculate and obtain the full symbol of the resolvent operator for the thermoelastic system {\rm\cite{LiuTan-22}} (it is also an elliptic system). Furthermore, the parabolic trace for this elliptic system can also be explicitly obtained.

    Generally, for any matrix-valued elliptic operator $L$, if one can explicitly obtain the inverse matrix for the principal symbol of the operator $(\tau I- L)^{-1}$, then all coefficients in the asymptotic expansion of the parabolic trace will  immediately be obtained by our new (algorithm) method.
\end{remark}

\begin{remark} Our paper {\rm\cite{Liu-21}} only studied the elastic Lam\'{e} system (i.e., elliptic system) for given $\mu>0$ and $\lambda+\mu>0$.
  As $\lambda+\mu\to 0$, our asymptotic expansion tends to the corresponding (well-known) result for the Laplace--Beltrami operator. In fact, all results and discussions in {\rm\cite{Liu-21}} still hold for $\mu>0$ and $\lambda+\mu\ge 0$.
\end{remark}

\begin{remark}
  Two-term asymptotics for the counting function $\mathcal{N}(\Lambda)$ as $\Lambda\to +\infty$ implies two-term asymptotics for the partition function $\mathcal{Z}(t)$ as $t\to 0^+$ under ``elastic billiards'' assumption, however, the converse may not true, where $\mathcal{N}(\Lambda):=\#\{k|\tau_k\le \Lambda\}$, and $\mathcal{Z}(t):=\sum_{k=1}^\infty e^{-t\tau_k}$. The reason is: the famous Tauberian theorem shows that one-term asymptotics for $\mathcal{Z}(t)$ as $t\to 0^+$ implies one-term asymptotics for $\mathcal{N}(\Lambda)$ as $\Lambda\to +\infty$. But, the Tauberian theorem is not true for two-term asymptotics. More precisely,
   let $v(s)$ vanish for $s<0$, be nondecreasing, continuous from the right, and suppose that for some constant $\alpha>0$,
   \begin{eqnarray*} \int_0^\infty e^{-ts} dv(s)\sim \frac{A}{t^\alpha} +\frac{B}{t^{\alpha-1}} +o(t^{1-\alpha})\;\;\;\,\mbox{as}\;\; t\to 0^+,\end{eqnarray*} then one can not obtain
   \begin{eqnarray*} v(s) \sim \frac{A}{\Gamma(\alpha +1)}s^\alpha +\frac{B}{\Gamma(\alpha)}s^{\alpha-1} +o(s^{\alpha-1})\;\;\;\mbox{as}\;\; s\to +\infty.\end{eqnarray*}
      Professor J. Korevaar (the author of ``Tauberian Theory: A Century of Developments'', Springer-Verlag, 2004) told the author of {\rm\cite{Liu-21}} this well-known fact with a detailed example in a private e-mail about sixteen years ago. Clearly, two-term asymptotics as $t\to 0^+$ for the partition function $\mathcal{Z}(t)$  in {\rm\cite{Liu-21}} can not go back to two-term asymptotics for $\mathcal{N}(\Lambda)$ as $\Lambda \to +\infty$. Thus, any numerical (verification) comparison must be done in the setting of  $\mathcal{Z}(t)$ as $t\to 0^+$ (here only an additional ``elastic billiards'' assumption is added) rather than to consider $\mathcal{N}(\Lambda)$ as $\Lambda\to +\infty$ as done in {\rm\cite{CaFrLeVa-22}}. Obviously, putting numerical experimental verification in the setting of $\mathcal{N} (\Lambda)$ with $\Lambda$ in a finite interval is a fundamental mistake in {\rm\cite{CaFrLeVa-22}}.

       In the numerical verification in {\rm\cite{CaFrLeVa-22}}, the number $\Lambda=1600$ or $\Lambda=2400$ or $\Lambda=2800$, or $\Lambda=2000$ or $\lambda=3000$ are too small. Even one takes $\Lambda=1000^{1000}$ or $\Lambda=10000^{10000}$ or $\Lambda=100000^{100000}$, it may not describe any true asymptotic case of $\mathcal{N}(\Lambda)-a\Lambda$  as $\Lambda\to +\infty$,  where $a$ is the one-term coefficient of $\mathcal{N}(\Lambda)$. Numerical experimental verification should consider the case after large $\Lambda$ rather than before. In {\rm\cite{CaFrLeVa-22}}, it seems that the numerical (experimental) data have ``satisfying'' effect for their second term of the counting function when $0\le \Lambda\le 3000$. However, it is useless to study the asymptotic behavior as $\Lambda\to +\infty$. Such a numerical experimental verification for small $\Lambda$ may not provide any valuable information for the future case of $\mathcal{N}(\Lambda)-a\Lambda$ as $\Lambda\to +\infty$. More precisely,  in order to verify $$\lim_{\Lambda\to +\infty} \frac{1}{\Lambda^{(n-1)/2} }\left[ \mathcal{N} (\Lambda)- a \mbox{Vol}_n (\Omega) \Lambda^{n/2} \right] = b'$$ for some constant $b'$, by the definition of limit, one should prove that for every number $\epsilon>0$ there is a number $M>0$ which depends on $\epsilon$ such that if $\Lambda>M$ then $$\bigg|\frac{1}{\Lambda^{(n-1)/2}}\left[ \mathcal{N} (\Lambda)- a \mbox{Vol}_n (\Omega) \Lambda^{n/2} \right]  - b'\bigg|<\epsilon.$$
       The authors of {\rm\cite{CaFrLeVa-22}} claimed that their result was ``correct'' only by ``numerically'' verifying the above inequality for some $\epsilon_0>0$ when $0\le \Lambda\le 3000$ or $0\le \Lambda\le 2400$ or $0\le \Lambda\le 2600$ for the unit disk and unit square. Such a so-called ``numerical verification'' method is not believable at all.
      \end{remark}

\begin{remark}  The authors of {\rm\cite{CaFrLeVa-22}} essentially discussed the elastic operator in a very special domain (i.e., the Euclidean upper half-space with flat metric $g_{jk}=\delta_{jk}$). Their method is ``stretch $\Omega$ by a linear factor $\kappa>0$, note that the
eigenvalues then rescale as $\kappa^{-2}$, and check the rescaling of the geometric invariants and of (1.20)'' (see, from line 1 to line 2 on p.$\,$6 of {\rm\cite{CaFrLeVa-22}}). However, for the two-term asymptotics, such a method will lead to losing much useful information because the boundary surface ``scrambles'' the elastic waves. In fact, in a Riemannian manifold $(\Omega, g)$, it follows from Lemma 2.1.1 of {\rm\cite{Liu-19}} that the elastic operator has the following local representation: \begin{eqnarray*} &&\mu\sum\limits_{j=1}^n\! \Big\{\! \Delta_g u^j \! +2 \!\sum\limits_{k,s,l=1}^n\! g^{kl} \Gamma_{sk}^j \frac{\partial u^s}{\partial x_l} \!+\!
   \sum_{k,s,l=1}^n \!\Big(g^{kl} \frac{\partial \Gamma^j_{sl}}{\partial x_k}\! +\!\sum\limits_{h=1}^n g^{kl} \Gamma_{hl}^j \Gamma_{sk}^h\! -\!\sum\limits_{h=1}^n g^{kl} \Gamma_{sh}^j \Gamma_{kl}^h \Big)u^s\Big\}\frac{\partial}{\partial x_j}\\
   &&\;\, +(\lambda+\mu) \operatorname{grad}\,\operatorname{div}\, \mathbf{u} \!+\! \mu\,\operatorname{Ric} (\mathbf{u}).\nonumber\end{eqnarray*}
In {\rm\cite{CaFrLeVa-22}}, by erroneously thinking of the above elastic operator defined on $(\Omega,g)$ as $\mu\sum_{j=1}^n \frac{\partial^2\mathbf{u}}{\partial x_j^2} +(\lambda+\mu) \nabla(\nabla\cdot \mathbf{u})$ on $\mathbb{\bar R}^{n}_{+}$, and further by mistaking a neighborhood of the boundary as the Euclidean upper half-space, the authors of {\rm\cite{CaFrLeVa-22}} had completely changed the original eigenvalues problems into other different spectral problems. The reader can not see where the branching Hamiltonian billiards condition (i.e., the corresponding billiard is neither dead-end nor absolutely periodic) on the Riemannian manifold is used in the proof of {\rm\cite{CaFrLeVa-22}}. In a Riemannian manifold, in order to get the two-term asymptotics of the counting function, the cut-off function (or pseudodifferential operator of order zero) should be introduced. See Response 6 in Section 3 for a correct approach to obtain the two-term asymptotics.
\end{remark}

\vskip 0.69 true cm

\section{Some responses}

\vskip 0.38 true cm

We will give some responses to remarks in \cite{CaFrLeVa-22} that involved the paper \cite{Liu-21}:

\vskip 0.25 true cm

\noindent {\bf Response 1}. \ {\it On p.$\,$27 of \cite{CaFrLeVa-22}, from line 6 to line 10, the authors of \cite{CaFrLeVa-22} wrote ``Finally, [Liu-21, p. 10166] contains a remarkable claim: ``McKean and Singer ...proved the famous Kac conjecture...''. Here the conjecture in question is the ``Can one hear the shape of a
drum?'' one, which was proven to be false in [GoWeWo-92]. The question of whether a similar conjecture is true for the operator of linear elasticity is still fully open, to the best of our knowledge.''}

\vskip 0.12 true cm
Clearly,  the authors of \cite{CaFrLeVa-22} have confused the famous Kac question and Kac conjecture, and given an incorrect comment.
Let us write out the original sentences in \cite{Kac-66} and \cite{MS-67} as follows.

   {\bf The Kac question}:  On p.\,2--p.\,3 of \cite{Kac-66}, M. Kac wrote ``Assume that for each $n$ the eigenvalue $\lambda_n$ for $\Omega_1$ is equal to the eigenvalue $\mu_n$ for $\Omega_2$.  Question: Are the regions $\Omega_1$ and $\Omega_2$ congruent in the sense of Euclidean geometry? ''

{\bf The Kac conjecture}:  On p.\,23 of \cite{Kac-66}, M. Kac wrote ``It is thus natural to conjecture that for a {\it smooth} drum with $r$ {\it smooth} holes
\begin{eqnarray*}  \sum_{n=1}^\infty e^{-\lambda_n t}\sim \frac{|\Omega|}{2\pi t}-\frac{L}{4} \frac{1}{\sqrt{2\pi t} }+\frac{1}{6}(1-r),
\end{eqnarray*}
and that therefore one can ``hear'' the connectivity of the drum!''

On p.\,44 of \cite{MS-67}, H. McKean and I. Singer wrote ``Kac was led to conjecture \begin{eqnarray*}(4b)\; \;\;\;\, \,\;\;\; \quad \;    Z=\frac{{\rm area}}{4\pi t}- \frac{{\rm length}/4}{\sqrt{4\pi t}}+\frac{1}{6} (1-h)+o(1) \;\;\;\;\;  (t\downarrow  0)
\end{eqnarray*}
for regions $D$ with smooth $B$ and $h < \infty$ holes, and was able to prove the correctness of the first $2$ terms. This jibes with an earlier conjecture of A. Pleijel and suggests that you can hear {\it the number of holes}. $(4b)$ will be proved below in a form applicable both to open manifolds with compact boundary and to closed manifolds.''

In \cite{Liu-21}, on p.$\,$10166, from line $\!-14$ to line $\!-1$, the author of
             \cite{Liu-21} accurately introduced the Kac question and the Kac conjecture according to the original literatures \cite{Kac-66} and \cite{MS-67}. The comment (in \cite{Liu-21}) ``.... a celebrated result of the spectral (geometric) invariants had been obtained by McKean and Singer. They proved the famous Kac conjecture and gave an explicit expression to the first three coefficients of asymptotic expansion for the heat trace of the Laplacian on a bounded domain $\Omega$ of a Riemannian manifold: ...'' is a  proper and correct statement.

\vskip 0.46 true cm

\noindent {\bf Response 2}. \    {\it In the footnote 3 on p.\,3 of {\cite{CaFrLeVa-22}}, the authors of {\cite{CaFrLeVa-22}} wrote ``We will not call (1.6) the Neumann condition in order to avoid confusion  with erroneous ``Neumann'' condition in \cite{Liu-21}.'' And for the remark on  p.\,9 of {\cite{CaFrLeVa-22}}, ``{\bf The first fundamental mistake} of {\cite{Liu-21}} is ... ''}

\vskip 0.20 true cm

In fact, in \cite{Liu-21} the Neumann boundary condition is $\frac{\partial \mathbf{u}}{\partial \boldsymbol{\nu}}:= 2\mu (\mbox{Def}\,\mathbf{u})^\#\boldsymbol{\nu} +\lambda(\mbox{div}\, \mathbf{u})\boldsymbol{\nu}$, as clearly pointed out from line 18 to line 21 on p.\,10166 of \cite{Liu-21}, and the author of \cite{Liu-21} wrote ``For the derivation of the Navier--Lam\'{e} elastic wave equations, its mechanical meaning, and the explanation of the Dirichlet and Neumann boundary conditions, we refer the reader to \cite{Liu-19} for the case of Riemannian manifold and ....''. In \cite{Liu-19}, the Neumann boundary condition is $2\mu (\mbox{Def}\, \mathbf{u})^\#\boldsymbol{\nu} + \lambda(\mbox{div}\, \mathbf{u})\boldsymbol{\nu}$ on $\partial \Omega$. \cite{Liu-19} is an earlier paper  which was posed on arXiv on Aug.\,14, 2019 by the author (Note that the paper \cite{Liu-21} was submitted to {\it The Journal Geometric Analysis} on Jul.\,21, 2020). In \cite{Liu-21},  the Neumann boundary condition is denoted by $\frac{\partial \mathbf{u}}{\partial \boldsymbol{\nu}}$. It is very clear from Lemma 2.1.1 of \cite{Liu-19} that $\frac{\partial \mathbf{u}}{\partial \boldsymbol{\nu}}:=2\mu (\mbox{Def}\, \mathbf{u})^\#\boldsymbol{\nu} +\lambda(\mbox{div}\, \mathbf{u})\boldsymbol{\nu}$, which equals to
$ \sum_{j,k=1}^n \big( \mu (u^j_{\;\,;k}+ u_k^{\;\,;j}) \nu^k +\lambda u^k_{\;\, ;k} \nu^j\big) \frac{\partial}{\partial x_j}$
(or equivalently, $\lambda  n^j \nabla_k u^k +\mu ( n^k \nabla_k u^j+ n_k \nabla_j u^k)$), where $\{\frac{\partial}{\partial x_j}\}_{j=1}^n$ is the nature coordinate basis. In addition, by applying the Green's formula to the elastic Lam\'{e} system, one can immediately see the corresponding Neumann boundary condition.

Let us point out the following facts (by applying Green's formula to the corresponding partial differential equation (or system of equations)):

\vskip 0.16 true cm

i).  \  For the Laplace equation $\Delta u=0$ in $\Omega$, the Dirichlet boundary condition is $u\big|_{\partial \Omega}$, and the Neumann boundary condition is $(\boldsymbol{\nu}\cdot \nabla u)\big|_{\partial \Omega}$, where $\boldsymbol{\nu}$ is the outward normal vector field on $\partial \Omega$;

\vskip 0.16 true cm

ii).  \  For the divergence type equation $\nabla \cdot (a(x) \nabla u)=0$ in $\Omega$, where $a(x)>0$ on $\bar \Omega$, the Dirichlet boundary condition is $u\big|_{\partial \Omega}$, and the Neumann boundary condition is $(\boldsymbol{\nu}\cdot a(x) \nabla u)\big|_{\partial \Omega}$;

\vskip 0.16 true cm

iii).  \  For the elastic Lam\'{e} system $\mu \Delta \mathbf{u} + (\mu +\lambda) \nabla (\nabla\cdot \mathbf{u}) =0$ in $\Omega$, where Lam\'{e} coefficients $\mu >0$ and $\mu+\lambda \ge 0$ on $\bar \Omega$,  the Dirichlet boundary condition is $\mathbf{u}\big|_{\partial \Omega}$, and the Neumann boundary condition is $\mu(\nabla \mathbf{u} +\nabla \mathbf{u}^T)^\#\boldsymbol{\nu} + \lambda(\mbox{div}\, \mathbf{u})\, \boldsymbol{\nu}$ on $\partial \Omega$ (see \cite{Liu-19});

\vskip 0.16 true cm

iv).  \  For the Maxwell's equations \begin{eqnarray*} \left\{ \begin{array}{ll} \mbox{curl} \mathbf{E} -i\omega \mu \mathbf{H}=0 \;\; &\mbox{in}\;\; \Omega,\\
\mbox{curl} \, \mathbf{H} +i \omega \sigma \mathbf{E} \;\; &\mbox{in}\;\; \Omega.\end{array}\right.\end{eqnarray*}
where $\omega>0$ is a fixed frequency; $\sigma$ is electronic coefficient and $\mu$ is the magnetic permeability,  the Dirichlet boundary condition is $ (\boldsymbol{\nu} \times \mathbf{E})\big|_{\partial \Omega}$, and the Neumann boundary condition is  $(\boldsymbol{\nu} \times \mathbf{H})\big|_{\partial \Omega}$ (see \cite{Liu-19.2});
\vskip 0.16 true cm
v).  \  For the Stokes flow equations (of a fluid) \begin{eqnarray*} \left\{ \begin{array}{ll} \mu \Delta \mathbf{u} +\nabla p=0 \;\; &\mbox{in}\;\; \Omega,\\
\mbox{div}\,  \mathbf{u} =0\;\; &\mbox{in}\;\; \Omega,\end{array}\right.\end{eqnarray*}
where $\mu >0$ is the viscosity constant,  the Dirichlet boundary condition is $\mathbf{u}\big|_{\partial \Omega}$, and the Neumann boundary condition is $\mu(\nabla \mathbf{u}+\nabla \mathbf{u}^T)^\# \boldsymbol{\nu} +p \boldsymbol{\nu}$ on $\partial \Omega$ (see \cite{Liu-20}).

\vskip 0.12 true cm
The boundary conditions depend strongly on the corresponding partial differential equation (or system of equations). They are an indivisible whole. Speaking boundary conditions without any partial differential equation(s) is not a proper statement. In \cite{Liu-21}, the author never said $\frac{\partial \mathbf{u}}{\partial \boldsymbol{\nu}}$ to be $\boldsymbol{\nu}\cdot \nabla \mathbf{u}$ as be understood in (1.33) on p.\,8 of \cite{CaFrLeVa-22}. Again, in \cite{Liu-21}, the Neumann boundary condition for elastic system is  $\mu(\nabla \mathbf{u}+\nabla \mathbf{u}^T)^\# \boldsymbol{\nu} +\lambda(\mbox{div}\,\mathbf{u}) \boldsymbol{\nu}$ as pointed out from line 18 to line 21 on p.\,10166 (see the paper \cite{Liu-19}, from version 1 to version 3).

On the one hand, in (1.1) and (1.7) of {\cite{CaFrLeVa-22}}, the authors of {\cite{CaFrLeVa-22}} essentially copied the expressions of \cite{Liu-19} for the Lam\'{e} operator and the Neumann boundary condition on a Riemannian manifold (see, Lemma 2.1.1 of \cite{Liu-19}) but did not cite (or mention) Lemma 2.1.1 of \cite{Liu-19}. On the other hand, in {\cite{CaFrLeVa-22}}, the authors of {\cite{CaFrLeVa-22}} (intentionally) wrongly explained the  Neumann boundary condition of \cite{Liu-21} and \cite{Liu-19}.

\vskip 0.49 true cm

\noindent {\bf Response 3}. \  {\it On p.\,26 of \cite{CaFrLeVa-22}, from line 20 to line 22, the authors of \cite{CaFrLeVa-22} wrote: ``Mathematically the reason for the erroneous expressions comes from the fact that Liu treats the principal symbol of $\mathfrak{L}$ as a {\it diagonal} matrix with eigenvalues (1.8) rather than a {\it diagonalisable} matrix with these eigenvalues.''}

This is not true. Clearly, the authors of \cite{CaFrLeVa-22} gave an irrelevant remark to the proof of \cite{Liu-21}.  Obviously, in \cite{Liu-21}, the principal symbol of $P_g$ (see, p.\,10177--p.\,10178 of \cite{Liu-21}) is
\begin{eqnarray*}
\mu\Big( \sum_{m,l=1}^n g^{ml} \xi_m\xi_l \Big){\mathbf{I}}_n +(\mu+\lambda)
\begin{bmatrix} \sum\limits_{m=1}^n g^{1m} \xi_m\xi_1 & \cdots & \sum\limits_{m=1}^n g^{1m} \xi_m\xi_n
  \\    \vdots & {} & \vdots  \\
 \sum\limits_{m=1}^n g^{nm} \xi_m\xi_1  & \cdots & \sum\limits_{m=1}^n g^{nm} \xi_m\xi_n  \end{bmatrix}, \end{eqnarray*}
 which is not a diagonal matrix. In \cite{Liu-21}, the author of \cite{Liu-21} never treated the principal symbol as being stated by the authors of \cite{CaFrLeVa-22}. Clearly, the authors of \cite{CaFrLeVa-22} misunderstood the main ideas of \cite{Liu-21}, and they tried to overthrow Theorem 1.1 of \cite{Liu-21} by their irrelevant remark.

\vskip 0.45 true cm

\noindent {\bf Response 4}. \   {\it  From line 3 to line 4 on p.\,27 of \cite{CaFrLeVa-22}, the authors of \cite{CaFrLeVa-22} wrote ``[Liu-21, formula (1.8)] uses the representation (1.24) in the Riemannian setting in an arbitrary
dimension. We recall that the definition of the operator curl is non-trivial in $d>3$ even in the
Euclidean case, and more clarification would be helpful here.''}
\vskip 0.12 true cm
In fact, from line 5 to line 17 on p.\,10166 of \cite{Liu-21}, the author of \cite{Liu-21} wrote ``In three spatial dimensions, .....''. It is clear that equation (1.8) in \cite{Liu-21} is considered in three dimensional case.

\vskip 0.49 true cm

\noindent {\bf Response 5}. \ {\it From line 23 to line 39 on p.$\,$26 of \cite{CaFrLeVa-22}, the authors of \cite{CaFrLeVa-22} wrote: ``We note that the expression for $\tilde{b}_{\mbox{Dir}}$  was already found in the 1960 paper by M. Dupuis, R.Mazo,
and L. Onsager \cite{DuMaOn-60}. Remarkably, this paper includes the critique of the 1950 paper by
E. W. Montroll who presented exactly Liu's expression (1.34) for the second asymptotic coefficient,
modulo some scaling, see [Mo-50, formulae (3)-(5)]. $\cdots$. In essence, Liu has rediscovered the seventy-year-old Montroll's formula which
had been known to be incorrect for sixty years.''}

\vskip 0.16 true cm

Obviously, the authors of \cite{CaFrLeVa-22} have misunderstood all results in these papers, and have  misled the reader by completely wrong statements. In \cite{Mo-50},
in order to discuss the effect of the volume and surface area to the heat capacities in low temperature,
E. W. Montroll gave the expression for the counting function of eigenvalues for the elastic normal modes of a three-dimensional rectangular solid $\{(x,y,z)\in \mathbb{R}^3\,\big|\, 0\le x\le L_x, 0\le y\le L_y, 0\le z\le L_z\}$ with the boundary condition $\boldsymbol{\nu} \cdot \nabla \mathbf{u}$ (according to the explanation of the boundary condition in \cite{CaFrLeVa-22}). In \cite{DuMaOn-60}, M. Dupuis, R. Mazo and L. Onsager investigated an isotropic solid at low temperatures whose model is a rectangular plate of thickness $l_3$ and other dimensions
$l_1$  and $l_2$ with realistic boundary conditions. The faces parallel to the plane of the plate are
supposed to be free of stresses, whereas the periodic
boundary conditions are given on the other faces (see, p.$\,$1453 of \cite{DuMaOn-60}).
  Clearly, these boundary conditions are completely different from the (whole) Dirichlet or Neumann (i.e., traction) boundary conditions in \cite{Liu-21}. In \cite{Mo-50} and \cite{DuMaOn-60}, the authors respectively calculated the counting functions with different boundary conditions for a very special domain (i.e., three-dimensional rectangular solid) by elementary calculations. But, in \cite{Liu-21}, we proved the asymptotic formulas for the heat traces of an elastic body  for general compact smooth manifold with smooth boundary (for the corresponding Dirichlet and the Neumann (i.e., traction) boundary conditions, respectively). The domains and the boundary conditions in three papers are quite different.

\vskip 0.39 true cm
\noindent {\bf Response 6}. \ {\it On p.\,27 of \cite{CaFrLeVa-22}, the authors of \cite{CaFrLeVa-22} wrote ``The paper does not refer to previous publications [De-12], [DuMaOn-60] or [SaVa-97].''}
\vskip 0.1 true cm
Clearly, the following two conjectures are still open:
 \vskip 0.1 true cm
 {\bf  Conjecture 1}:  \ Let $\Omega$ be a compact, connected smooth Riemannian $n$-manifold with smooth boundary $\partial \Omega$.
 The following two-term asymptotics hold:
 \begin{eqnarray} \label{2022.9.7-1} &&\;\;\, \, \mathcal{N}_{\mp} (\Lambda) \sim \frac{1}{ \Gamma(1+\frac{n}{2})} \Big( \frac{n-1}{(4\pi \mu)^{n/2}} +\frac{1}{(4\pi (\lambda+2\mu))^{n/2}}\Big)\mbox{vol}_n(\Omega) \Lambda^{n/2}\\
&& \quad \;  \mp \frac{1}{4 \,\Gamma(1+\frac{n-1}{2})} \Big( \frac{n-1}{(4\pi\mu)^{(n-1)/2}} +\frac{1}{(4\pi(\lambda+2\mu))^{(n-1)/2}}\Big)\mbox{vol}_{n-1}(\partial \Omega) \,\Lambda^{(n-1)/2}\; \quad \;\;  \mbox{as}\;\, \Lambda \to +\infty, \nonumber \end{eqnarray}
 where $\mathcal{N}_{-} (\Lambda):=\#\{ k\big| \tau^{-}_k \le \Lambda\}$ (respectively, $\mathcal{N}_{+} (\Lambda):=\#\{ k\big| \tau_k^{+}\le \Lambda\}$) is the counting function of elastic eigenvalues for the Dirichlet (respectively, Neumann) boundary condition.

\vskip 0.1 true cm
 When $\lambda+\mu=0$, Conjecture 1 just is the famous Weyl conjecture, which has about 110-year history and remains open.

  \vskip 0.2 true cm

 {\bf Conjecture 2}:  Let $\Omega$ be a compact, connected smooth Riemannian $n$-manifold with smooth boundary $\partial \Omega$. Suppose that the corresponding elastic billiards is neither dead-end nor absolutely periodic for the $\Omega$. The above two-term asymptotics (\ref{2022.9.7-1}) hold.

 \vskip 0.2 true cm

 Obviously,  for any  given compact smooth Riemannian manifold $\Omega$ with smooth boundary, under the elastic billiards condition, an explicit two-term  asymptotic expansion of the counting function $\mathcal{N}_{\mp}(\Lambda)$ as $\Lambda\to +\infty$  had not been given in  \cite{SaVa-97}. In [De-12] and  [DuMaOn-60], the counting functions of elastic eigenvalues were explicitly calculated only for very special domain (a three-dimensional rectangular solid) with other boundary conditions.
 However, in Theorem 1.1 of \cite{Liu-21}, our two-term asymptoic formulas of the heat trace hold for any given compact connected smooth Riemannian manifold $\Omega$ with smooth boundary (do not need any additional assumption).

  Similar to Ivrii's method in \cite{Iv-80} (or \cite{Iv-16}) and H\"{o}rmander's technique in \cite{Ho-85}, it is not difficult to prove the above Conjecture 2. In this new approach which is completely different from \cite{SaVa-97} and \cite{CaFrLeVa-22}, the spectral measure and its Fourier transform will be introduced; the fundamental solution $\mathbf{u}(t,x,y)$ of the elastic wave equations with suitable boundary conditions will be considered; the asymptotic expansion of the normal singularity of $\sigma (t)=\int_\Omega \mbox{Tr}\,(\mathbf{u} (t,x,x))\,dx$ as well as $\sigma_\zeta (t)= \int_\Omega \mbox{Tr}\,(\zeta (x) \mathbf{u}(t,x,x)) dx$ will be estimated and explicitly calculated; and the branching Hamiltonian billiards condition on Riemannian manifold will be well used.

     However, not every manifold satisfies the billiards condition. For example, when $\lambda+\mu=0$ (the corresponding elastic operator becomes the Laplacian), the semi-sphere $\mathbb{S}^{n-1}_{+}:=\{(x_1,\cdots, x_{n-1}$, $x_n)\in \mathbb{R}^n \big|x_1^2 +\cdots +x_{n-1}^2+x_n^2=1,\; x_n\ge 0\}$ does not satisfy the billiards condition because every point $(x,\xi)\in S^* (\mathbb{S}^{n-1}_{+})$  is periodic  (Here $S^* (\mathbb{S}^{n-1}_{+})$ is the fiber of cotangent unit sphere over $\mathbb{S}^{n-1}_{+}$).
\vskip 0.1 true cm

\addvspace{16mm}

\section*{Acknowledgements}

\addvspace{6mm}

This research was supported by NNSF of China (11671033/A010802) and NNSF of China \\
(11171023/A010801).

\addvspace{12mm}

\addvspace{5mm}

\end{document}